\documentclass[10pt,fleqn]{amsart}

\usepackage{amsmath,amssymb,latexsym}
\usepackage[mathscr]{eucal}

\theoremstyle{plain}
\newtheorem{theorem}{Theorem}
\newtheorem{lemma}{Lemma}

\numberwithin{equation}{section}

\newcommand{\dy}{\partial}

\newcommand{\ddt}[1]{\frac{{\rm d}{#1}}{{\rm d}{t}}}

\newcommand{\Zahl}{\mathbb{Z}}

\newcommand{\Comp}{\mathbb{C}}

\newcommand{\e}{{\rm e}}
\newcommand{\Chi}{{\rm X}}
\newcommand{\eps}{\varepsilon}

\newcommand{\lapl}[1]{\Delta_{#1}^{}}
\newcommand{\ilapl}[1]{\Delta_{#1}^{-1}}

\newcommand{\ub}{\boldsymbol{u}}
\newcommand{\vb}{\boldsymbol{v}}
\newcommand{\vbb}{\bar{\vb}}
\newcommand{\Vbb}{{\bar{\boldsymbol{V}}}}
\newcommand{\xb}{\boldsymbol{x}}
\newcommand{\kb}{{\boldsymbol{k}}}
\newcommand{\gb}{\boldsymbol{\nabla}}
\newcommand{\sgb}{\boldsymbol{\nabla}^\perp}

\newcommand{\sle}{{\scriptscriptstyle<}}

\newcommand{\Pl}{{\sf P}^{\!{}^<}}
\newcommand{\Pz}{{\sf P}_{\!z}}

\newcommand{\qt}{\tilde q}

\newcommand{\Gy}{\mathcal{G}}
\newcommand{\Gyle}{\Gy^{\!{}^<}}
\newcommand{\Gylh}{{\hat\Gy}}
\newcommand{\Fy}{\mathcal{F}}
\newcommand{\Sy}{\mathcal{S}}

\newcommand{\Dom}{\mathscr{M}}

\newcommand{\qa}{\tilde q}
\newcommand{\xa}{\chi}
\newcommand{\fa}{\phi}

\newcommand{\fv}{f_{\vb}^{}}
\newcommand{\fr}{f_\rho^{}}
\newcommand{\fq}{f_q^{}}
\newcommand{\fh}{f_\chi^{}}
\newcommand{\ff}{f_\phi^{}}

\newcommand{\uba}{\ub}
\newcommand{\vba}{\vb}
\newcommand{\rba}{\rho}

\newcommand{\cnst}[1]{c_{#1}^{}}

\newcommand{\Rr}{\mathcal{R}}
\newcommand{\Rv}{\mathcal{R}_{\bar{\vb}}}
\newcommand{\Rx}{\mathcal{R}_\chi}
\newcommand{\Rf}{\mathcal{R}_\phi}

\newcommand{\FD}{\mathsf{D}}

\newcommand{\dl}{\delta}

\begin{document}

\title[Exponential Approximations for Primitive Equations]%
{Exponential Approximations for the\\
Primitive Equations of the Ocean}

\author[Temam]{R.~Temam}
\email{temam@indiana.edu}
\urladdr{http://mypage.iu.edu/\~{}temam}

\author[Wirosoetisno]{D.~Wirosoetisno}
\email{djoko.wirosoetisno@durham.ac.uk}
\urladdr{http://www.maths.dur.ac.uk/\~{}dma0dw}
\curraddr[Wirosoetisno]{Department of Mathematical Sciences\\
   University of Durham\\
   Durham, DH1~3LE, United Kingdom}

\address{The Institute for Scientific Computing and Applied Mathematics\\
   Indiana University, Rawles Hall\\
   Bloomington, IN~47405--7106, United States}

\thanks{This research was supported by grants NSF~0305110,
DOE~DE--FG02--01ER63251:A000 and the Research Fund of Indiana University}

\keywords{Singular perturbation, exponential asymptotics, Gevrey regularity, primitive equations}
\subjclass[2000]{Primary: 35B25, 76U05}


\begin{abstract}
We show that in the limit of small Rossby number $\eps$, the primitive
equations of the ocean (OPEs) can be approximated by ``higher-order
quasi-geostrophic equations'' up to an exponential accuracy in $\eps$.
This approximation assumes well-prepared initial data and is valid for
a timescale of order one (independent of $\eps$).
Our construction uses Gevrey regularity of the OPEs and a classical
method to bound errors in higher-order perturbation theory.
\end{abstract}

\maketitle


\section{Introduction}\label{s:intro}

We consider the primitive equations for the ocean (henceforth OPEs),
scaled as in \cite{petcu-temam-dw:rgpe}
\begin{equation}\label{q:pe-uvr}\begin{aligned}
   &\dy_t\vb + \frac1\eps \bigl[ \vb^\perp + \gb p \bigr]
	+ \ub\cdot\gb \vb = \mu\lapl3 \vb + \fv,\\
   &\dy_t\rho - \frac1\eps w + \ub\cdot\gb \rho = \mu\lapl3 \rho + \fr,\\
   &\gb\cdot\ub = \gb\cdot\vb + w_z = 0,\\
   &\rho = -p_z.
\end{aligned}\end{equation}
Here $\ub=(u,v,w)$ is the three-dimensional fluid velocity, with $\vb=(u,v)$
its horizontal component and $\vb^\perp=(-v,u)$;
$p$ is the pressure;
$\rho$ is the perturbation density (not including the mean stable
stratification which figures into the $\eps$ in the equation for $\rho$).
We write $\gb_2:=(\dy_x,\dy_y)$ and $\gb_3:=(\dy_x,\dy_y,\dy_z)$;
when no ambiguity may arise, we simply write $\gb$.
Similarly, we write $\lapl2:=\dy_x^2+\dy_y^2$ and
$\lapl3:=\dy_x^2+\dy_y^2+\dy_z^2$.
The parameter $\eps$ is related to the Rossby and Froude numbers;
in this article we shall be concerned with the limit $\eps\to0$,
and for convenience we assume that $\eps\le1$ (further restrictions
on $\eps$ will be stated below).
In general the viscosity coefficients for $\vb$ and $\rho$ are different;
we have set them all to $\mu$ for clarity of presentation (the general
case does not introduce any more essential difficulty).
The forcings $\fv$ and $\fr$ are assumed to be independent of time.

We work in three spatial dimensions,
$\xb = (x,y,z) \in [0,L_1]\times[0,L_2]\times[-L_3/2,L_3/2]$\penalty0$=: \Dom$,
with periodic boundary conditions assumed.
Following common practice in numerical simulations of stratified turbulence
(see, e.g., \cite{bartello:95}), the dependent variables are assumed to have
the following symmetries:
\begin{equation}\label{q:pe-sym}\begin{aligned}
   &\vb(x,y,-z) = \vb(x,y,z),
	&\qquad
   &p(x,y,-z) = p(x,y,z),\\
   &w(x,y,-z) = -w(x,y,z),
	&\qquad
   &\rho(x,y,-z) = -\rho(x,y,z);
\end{aligned}\end{equation}
we say that $\vb$ and $p$ are {\em even} in $z$,
while $w$ and $\rho$ are {\em odd} in $z$.
If in addition $\fv$ is even and $\fr$ is odd in $z$, it can be verified
that this symmetry is preserved by the OPEs \eqref{q:pe-uvr}, that is,
if it holds at $t=0$, it continues to hold for $t>0$.
Since $w$ and $p$ are periodic in $z$, we have $w(x,y,-L_3/2)=w(x,y,L_3/2)=0$
and $\rho(x,y,-L_3/2)=\rho(x,y,L_3/2)=0$;
similarly, $u_z^{}=0$, $v_z^{}=0$ and $p_z^{}=0$ on $z=0,\pm L_3/2$
if they are sufficiently smooth (as will be assumed below).
One may consider the symmetry conditions \eqref{q:pe-sym} as a way to
impose the boundary conditions $w=0$, $\rho=0$, $u_z^{}=0$, $v_z=0$ and
$p_z^{}=0$ on both $z=0$ and $z=L_3/2$, in the effective domain
$[0,L_1]\times[0,L_2]\times[0,L_3/2]$.
All variables and the forcing are assumed to have zero mean over $\Dom$;
the symmetry conditions above ensure that this also holds for their products
that appear below.


It is known that, given sufficiently regular initial data, the OPEs have a
unique strong solution for all time \cite{cao-titi:u-3dpe},\cite{kobelkov:06};
the results of \cite{petcu-dw:gev3} then imply that a very regular
solution, belonging to a Gevrey space defined below, exists for all time.
The existence of finite-dimensional global attractors in various spaces
(see \cite{temam:iddsmp}), and thus uniform boundedness of the solution,
have been shown in \cite{ju:3dpe};
see also the result in \cite{petcu:3dpe}.
One can thus regard the OPEs as equivalent to a finite-dimensional system.
This fact is an important ingredient for our main result.


\medskip
As noted earlier, we are concerned with the limit $\eps\to0$.
A flow is said to be in {\em geostrophic balance\/} if the $\mathcal{O}(1/\eps)$
terms in (\ref{q:pe-uvr}a) vanish, i.e.\ if
\begin{equation}\label{q:geobal}
   \vb = \sgb p
\end{equation}
where $\sgb:=(-\dy_y,\dy_x)$.
Note that this implies $\sgb.\vb=0$, which with the fact that $w$ is
odd in $z$ in turn implies that $w=0$, so the $\mathcal{O}(1/\eps)$ term
in (\ref{q:pe-uvr}b) also vanishes.
Assuming that \eqref{q:geobal} is satisfied initially, it is well known
that the OPEs can be approximated by the simpler quasi-geostrophic
equation (QGE),
\begin{equation}\label{q:qge}
   \dy_t q^{\rm g} + \vb^{\rm g}\cdot\gb q^{\rm g} = \mu \lapl3 q^{\rm g} + \fq\,,
\end{equation}
which only involves a single variable $q^{\rm g}$;
here $\fq:=\sgb\cdot\fv-\dy_z\fr$.
The original variables are recovered using
\begin{equation}\label{q:qg-inv}
    \vb^{\rm g}=\sgb\ilapl3 q^{\rm g}
    \quad\textrm{and}\quad
    \rho^{\rm g}=-\dy_z\ilapl3 q^{\rm g},
\end{equation}
Here $\ilapl3$ is uniquely defined to have zero mean over $\Dom$.
The quantities $\vb^{\rm g}$ and $\rho^{\rm g}$ derived from the potential
vorticity $q^{\rm g}$ are said to be {\em geostrophic\/}.
The fact that $\dy_z w+\gb\cdot\vb=0$ and (\ref{q:qg-inv}a) imply
that $w^{\rm g}=0$.

An approximation result was obtained in \cite{bourgeois-beale:94}
for a closely related system where it was shown that the QGE
\eqref{q:qge} with $\mu=0$ and $f_q=0$ approximates (the unforced,
inviscid and Boussinesq analogue of) the OPE \eqref{q:pe-uvr} with
an error of order $\eps$ for $t\in[0,T]$, viz.,
\begin{equation}\label{q:bourgeois-beale}
   |\vb(t)-\sgb\ilapl3 q^{\rm g}(t)|_{H^4(\Dom)} + |w(t)|_{H^4(\Dom)}
   + |\rho(t)+\dy_z\ilapl3 q^{\rm g}(t)|_{H^4(\Dom)} \le \eps
\end{equation}
assuming that the left-hand side less than $\eps$ at $t=0$ (and given
sufficient regularity).
We shall not attempt to follow their approach in this article as it
makes no qualitative difference for our exponential-order result.

Regarding this result as a first-order approximation in the parameter $\eps$,
a natural question is whether one can obtain higher-order approximations.
Put differently, we would like to know how fast, as $\eps\to0$, the solution
of the PE converges to the solution of a simpler system analogous to the
QGE \eqref{q:qge}.
Our main purpose in this article is to show that convergence of any algebraic
order is possible for $\eps$ sufficiently small, resulting in an error
estimate which is exponentially small in $\eps$.

\medskip
In the geophysical parlance, the dynamics of a fluid flow is said to be
{\em balanced\/} if the solution stays near a subspace (``balance manifold'')
for some time.
The {\em order\/} of the balance dynamics measures how close this
approximation is in terms of $\eps$.
A result such as \eqref{q:bourgeois-beale} implies that the
quasi-geostrophic relations \eqref{q:qg-inv} define a balance
dynamics of order $\eps$.
Our result (Theorem~\ref{t:main} below) thus implies that a {\em balanced
dynamics of exponential order\/} exists for a timescale of order one
(i.e.\ independent of the Rossby number).
Note, however, that our result is actually stronger than the geophysical
definition of balance, since it gives a pointwise approximation
to an exponential order.

Also of geophysical interest is whether the solution (when appropriately
initialised) stays close to a balance manifold over longer timescales.
Using careful estimates (of a different type than those used here), it may
be possible to prove that the solution remains ``exponentially balanced''
for timescales of order $1/\eps$ (although it is clear that one cannot
expect pointwise accuracy over this timescale);
we plan to report on this in a future work \cite{temam-dw:lbal}.


\medskip
We begin with some notations.
We write the horizontal velocity $\vb$ as
\begin{equation}\label{q:hodge}
   \vb(x,y,z) = \vbb(z) + \sgb\psi + \gb\chi
\end{equation}
where $\vbb(z)$ is the mean vertical shear, i.e.\ the $(x,y)$-average
of $\vb$.
We introduce the {\em linearised potential vorticity\/} (which has appeared
in the QGE above)
\begin{equation}
   q := \sgb\cdot\vb-\rho_z \,.
\end{equation}
The streamfunction $\psi$ and velocity potential $\chi$ are defined
as follows.
First, let
\begin{equation}\label{q:chi-def}
   \chi := \ilapl2\gb\cdot\vb,
\end{equation}
where, here and henceforth, $\ilapl2$ is uniquely defined to have zero
mean on each $(x,y)$-plane.
Next, let
\begin{equation}\label{q:ur-qxf}
   \phi(x,y,z) := \int_0^z \ilapl2 \bigl[ \lapl3\rho + \dy_{z'} q \bigr]
	\;{\rm d}z'.
\end{equation}
And finally, let
\begin{equation}\label{q:psi-def}
   \psi := \ilapl3\bigl[ q + \phi_{zz} \bigr],
\end{equation}
where $\ilapl3$ is uniquely defined to have zero average over $\Dom$.
This completes the definition of $(q,\vbb,\chi,\phi)$ in terms of
$(\vb,\rho)$.

It can be verified that these definitions imply the following relations,
\begin{equation}\label{q:misc}
   \lapl2\psi = \sgb\cdot\vb,\qquad
   p = \psi - \phi
   \qquad\textrm{and}\qquad
   \lapl2\chi = \gb\cdot\vb.
\end{equation}
The second equation is useful to compute $\rho=-\dy_z p$ from $q$ and $\phi$;
$\vb$ can be computed from $(q,\vbb,\chi,\phi)$ using \eqref{q:hodge} and
\eqref{q:psi-def}.
Using the incompressibility condition $\gb\cdot\vb+w_z=0$ and the fact
that $w(x,y,0)=0$, we also have
\begin{equation}\label{q:wdef}
   w(x,y,z) = - \int_0^z \lapl2\chi(x,y,z') \;{\rm d}z'.
\end{equation}

The primitive equations \eqref{q:pe-uvr} can now be written in
the following form,
{\renewcommand\minalignsep{0pt}
\begin{equation}\label{q:pe-qxf}\begin{aligned}
   &\dy_t q && &&+ \sgb\cdot(\ub\cdot\gb\vb) - \dy_z(\ub\cdot\gb\rho)
	&&= \mu\lapl3 q + \fq,\\
   &\dy_t \vbb_z &&+\frac1\eps \vbb_z^\perp
	&&+ \dy_z\bigl(\overline{\vphantom{+}w\vb_z}\bigr)
        &&= \mu \vbb_{zzz} + \dy_z\bar\fv,\\
   &\dy_t\lapl3\chi &&-\frac1\eps\lapl3\phi
        &&+ \lapl3\ilapl2\gb\cdot(\ub\cdot\gb\vb)
	&&= \mu\Delta_3^2\chi + \lapl3 \fh,\\
   &\dy_t \phi_{zz} &&+ \frac1\eps \lapl3\chi
        &&+ \dy_{zz}\ilapl2 \sgb\cdot(\ub\cdot\gb \vb)
	+ \rlap{$\Pz\dy_z(\ub\cdot\gb\rho)$}\\
	& && && &&= \mu\lapl3\phi_{zz} + \dy_{zz}\ff,
\end{aligned}\end{equation}}%
where overbar denotes average over $(x,y)$ and $\Pz$ is a projection which
removes the $(x,y)$-average of a function, $\Pz\rho:=\rho-\bar\rho$, etc.
The forcing terms in \eqref{q:pe-qxf} are related to those in \eqref{q:pe-uvr}
by
\begin{equation}\label{q:fuvr-fqxf}\begin{aligned}
   &\fq = \sgb\cdot\fv - \dy_z\fr\,,
   &\qquad &f_\chi = \ilapl2\gb\cdot\fv\,,\\
   &\bar\fv = \overline{\fv}\,,
   &&\dy_{zz} f_\phi = \ilapl2\dy_{zz}\sgb\cdot\fv + \Pz\dy_z\fr.
\end{aligned}\end{equation}

It is easily verified that the $\mathcal{O}(1/\eps)$ terms in
\eqref{q:pe-qxf} are antisymmetric.
If the fast variables $(\vbb,\chi,\phi)$ are all zero---a condition that
is equivalent to the geostrophic balance relation \eqref{q:geobal}---they
can only grow through the nonlinear terms involving the slow variable $q$
and the nonlinear terms on the right-hand side.
Further analysis (cf.\ \cite{bourgeois-beale:94},\cite{temam-dw:lbal})
shows that they remain small over a timescale of order one.
Our main result [Theorem~\ref{t:main}] states that, for certain initial
data ($\vbb$, $\chi$ and $\phi$ given in terms of $q$), this system can
be approximated up to an error exponentially small in $\eps$ by the
solution of an equation for a finite-dimensional slow variable $\qt$.

We note that the form \eqref{q:pe-qxf} will only be used for (i)~the
construction of the approximation $\vb^*(\qt;\eps)$ and $\rho^*(\qt;\eps)$,
and (ii)~the integration of the finite-dimensional variable $\qt$.
The approximation error estimates are in terms of $(\vb,\rho)$ and
$(\vb^*,\rho^*)$, so we will not need a priori estimates for \eqref{q:pe-qxf}.


\section{Gevrey Regularity and Statement of the Main Result}\label{s:gevrey}

Let $W=(u,v,\rho)$.
The periodic boundary conditions allow us to write functions (assuming
sufficient regularity) in Fourier series,
\begin{equation}
   W(\xb,t) = \sum_{\kb\in\Zahl^3} W_\kb(t)\, {\rm e}^{{\rm i}\kb'\cdot\xb},
\end{equation}
where $\kb':=(k_1',k_2',k_3')$ with $k_i':=2\pi\,k_i/L_i$.
For $\sigma>0$ fixed, the Gevrey norm $\|\cdot\|_\sigma$ is defined as
\begin{equation}\label{q:gevdef}
   \|W\|_\sigma^2 := \sum_{\kb\in\Zahl^3} {\rm e}^{2\sigma|\kb|} |W_\kb|^2.
\end{equation}
It follows from this that if $\|W\|_\sigma<\infty$, $|W_\kb|$
decays exponentially fast in $|\kb|$;
this fact will play an important role below.

The following result on the OPEs is proved in \cite{petcu-dw:gev3}
(this has been extended to all time in \cite{petcu:3dpe}, but the
version here is more useful for our present need):

\begin{lemma}\label{t:gev}
Let $W_0=(u_0,v_0,\rho_0)$ be such that $|\gb_3 W_0| = C_{1,0} < \infty$
and suppose that $|\gb_3\fv|+|\gb_3\fr| = C_f <\infty$.
Then for any $C_1>C_{1,0}$ there exists a
$T_1(C_{1,0},C_f,C_1;\sigma,\mu,\Dom)>0$ such that for $t\in[0,T_1]$,
\begin{equation}\label{q:bdW1}
   |\gb_3 W(t)| \le C_1.
\end{equation}
Now let $\sigma>0$ be fixed, and let $W_0$ be such
that $\|\lapl3 W_0\|_\sigma = C_{\sigma,0} <\infty$;
in addition, suppose that
$\|\gb_3\fv\|_\sigma+\|\gb_3\fr\|_\sigma = C_f' <\infty$.
Then for any $C_\sigma>C_{\sigma,0}$, there exists a
$T_\sigma(C_{\sigma,0},C_f',C_\sigma;\sigma,\mu,\Dom)>0$ such that for
$t\in[0,T_\sigma]$,
\begin{equation}\label{q:bdWsig}
   \|\lapl3 W(t)\|_\sigma^{} \le C_\sigma.
\end{equation}
\end{lemma}

\noindent
We stress that $T_1$ and $T_\sigma$ are independent of $\eps$.
(Here and henceforth, all constants are understood to be positive.)

Given a fixed $\kappa>0$, we define a finite-dimensional truncation
of $W$ as
\begin{equation}
   W^<(\xb,t) := (\Pl W)(\xb,t)
	:= \sum_{|\kb|<\kappa} W_\kb(t) {\rm e}^{{\rm i}\kb'\cdot\xb},
\end{equation}
and $W^> := W - W^<$.
It is clear from this definition that the projection $\Pl$ is orthogonal.
For the low modes $W^<$ we have the ``reverse Poincar{\'e} inequality'',
\begin{equation}\label{q:ipoi}
   |\gb_3 W^<|_s \le c\,\kappa\, |W^<|_s.
\end{equation}
Here and throughout this article, $\|\cdot\|_\sigma$ denotes the Gevrey
norm \eqref{q:gevdef} and $|\cdot|_s$ denotes the usual Sobolev $H^s$ norm;
when no confusion may arise, we will often write simply $|\cdot|$ for the
$L^2$ norm $|\cdot|_0^{}$.
If $\|W\|_\sigma<\infty$, the exponential decay of the Fourier coefficients
$W_\kb$ allows us to bound the high modes $W^>$ in $L^2$ as
\begin{equation}\label{q:whisig}
   |W^>|_0 \le C_0\, {\rm e}^{-\sigma\kappa} \|W\|_\sigma
\end{equation}
or, more generally in $H^s$,
\begin{equation}\label{q:whisigs}
   |W^>|_s^{} \le C_s \kappa^s{\rm e}^{-\sigma\kappa} \|W\|_\sigma.
\end{equation}

We can now state our main result:

\begin{theorem}\label{t:main}
Let $\sigma\in(0,1)$ be fixed.
Suppose that the initial data $W_0=(u_0,v_0,\rho_0)$ satisfies
\begin{equation}
   \|\gb_3^3 W_0\|_\sigma^{} < \infty,
\end{equation}
and the forcing satisfies
\begin{equation}
   \|\gb_3^3\fv\|_\sigma + \|\gb_3^3\fr\|_\sigma<\infty\,.
\end{equation}
Then there exists an
$\eps_0(\|\gb_3^3 W_0\|_\sigma^{},\|\gb_3^3f\|_\sigma^{};\sigma,\mu,\Dom)>0$,
and a $\qt$ (which is finite dimensional) evolving according to
\begin{equation}
   \dy_t \qt + \Pl \bigl[ \sgb\cdot(\ub^*\cdot\gb\vb^*)
	- \dy_z(\ub^*\cdot\gb\rho^*) \bigr]
	= \mu\lapl3\qt + f_q^<\,,
\end{equation}
where $\ub^*=\ub^*(\qa;\eps)$ and $\rho^*=\rho^*(\qa;\eps)$ are constructed
in the proof below, such that if $\eps\in(0,\eps_0)$ and the initial data
satisfies
\begin{equation}\label{q:idhypo}
   |\vb_0-\vb^*(\qt(0);\eps)|_0^2 + |\rho_0-\rho^*(\qt(0);\eps)|_0^2
	\le C_{\rm id} \exp\bigl(-2\sigma/\eps^{1/4}),
\end{equation}
for some constant $C_{\rm id}$ {\rm (\/}independent of $\eps$
and $W_0${\/\rm)\/} then a similar estimate,
\begin{equation}\label{q:mainbd}
   |\vb(t)-\vb^*(\qt(t);\eps)|_0^2 + |\rho(t)-\rho^*(\qt(t);\eps)|_0^2
	\le 4C_{\rm id} \exp\bigl(-2\sigma/\eps^{1/4}),
\end{equation}
holds for $t\in[0,T_*]$ where
$T_*=T_*(\|\gb_3^3 W_0\|_\sigma^{},\|\gb_3^3f\|_\sigma^{};\sigma,\mu,\Dom)$.
\end{theorem}

\medskip\noindent{\bf Remarks.} 

\noindent{\bf 1.} The $\vb^*$ and $\rho^*$ are the ``exponential order''
analogues of the (leading-order) quasi-geostrophic $\vb^{\rm g}$ and
$\rho^{\rm g}$.
Unlike the latter, which is infinite-dimensional, our higher-order
approximations are based on a finite-dimensional variable and depend on
Gevrey regularity of the parent system \eqref{q:pe-uvr}.

\noindent{\bf 2.} The sharpness of our bound \eqref{q:mainbd} is unclear,
but an explicit example constructed in \cite{jv-yavneh:04} for a closely
related problem suggests that one cannot do better than $\exp(-c/\eps)$.

\noindent{\bf 3.} Although in this article we only treat the primitive
equations for the ocean \eqref{q:pe-uvr}, the method described here
should be applicable to many PDEs with a small parameter for which a
Gevrey regularity result can be proved.


\section{Proof of Theorem~\ref{t:main}}\label{s:pf}

Our approach is inspired by \cite{matthies:01}, which is an averaging-type
result on {\em single-frequency\/} infinite-dimensional systems,
and by \cite[App.~B]{cotter:th}, which considers a finite-dimensional problem.
Our present problem is infinite-dimensional with an infinite number of
frequencies, so small denominator problems would appear in an extension
of \cite{matthies:01} to this case.
However, by considering the singular perturbation problem (hence the
requirement for special initial data) in this article, we can avoid the small
denominator problem.

\medskip
\noindent The plan of the proof is as follows:

\medskip\noindent * We will build up our $n$th-order approximation $W^n$ as:
\begin{align}
   &W = W^< + W^>\\
   &W^< = W^n + \hat W
\end{align}
The high modes $W^>$ are exponentially small thanks to the Gevrey
regularity of $W$;
\eqref{q:whisigs} then implies that they are exponentially small in
any Sobolev norm.
The low modes $W^<$ (which are finite-dimensional) can then be approximated
by $W^n$, which is constructed such that $\hat W$ is exponentially
small in $\eps$.
The total error $W-W^n=W^>+\hat W$ is then also exponentially small.
Like $W^<$, the approximate solution $W^n=(\vb^n,\rho^n)$ is
finite dimensional,
\begin{equation}
   W^n(\xb,t) = \sum_{|\kb|<\kappa}
	W^n_\kb(t) \e^{{\rm i}\kb\cdot\xb}.
\end{equation}

\medskip\noindent * It is ``slaved'' to the slow variable $\qa$
in the following manner:
\begin{align}
   &\vba^n = \Vbb^n + \sgb\psi^n + \gb\Chi^n
	\quad\textrm{with}\quad
	\lapl3\psi^n = \qa + \dy_{zz}\Phi^n \label{q:vbn}\\
   &\rho^n = \ilapl3\bigl[ \lapl2\dy_z\Phi^n - \dy_z\qa \bigr],\label{q:rn}
\end{align}
where $\Vbb(\qa;\eps)$, $\Chi^n(\qa;\eps)$ and $\Phi^n(\qa;\eps)$ are
functions to be computed below.
For a motivation of this ``slaving ansatz'', we refer the reader
to \cite{dw-tgs-temam:wwrg}.
Even when working with finite-dimensional systems, we shall keep the
convenient PDE notation, remembering that in this case we have the
reverse Poincar{\'e} inequality \eqref{q:ipoi}.
We note for future reference
\begin{equation}\label{q:wn}
   w^n = -{\int_0^z} \lapl2\Chi^n(\qa;\eps) \>{\rm d}z'.
\end{equation}

\medskip\noindent * The (truncated) potential vorticity $\qa$ evolves
according to our ``higher-order quasi-geostrophic equations'' (QGE${}^n$),
\begin{equation}\label{q:qgen}
   \dy_t\qa + \Pl\bigl[ \sgb\cdot(\ub^n\cdot\gb\vb^n)
	- \dy_z(\ub^n\cdot\gb\rho^n) \bigr]
	= \mu \lapl3\qa + f_q^< \,,
\end{equation}
where $\ub^n=(\vb^n,w^n)$ and $\rho^n$ are determined in terms of
$(\qa,\Vbb^n,\Chi^n,\Phi^n)$ by \eqref{q:vbn}--\eqref{q:wn}.

\medskip\noindent * In turn, $(\Vbb^n,\Chi^n,\Phi^n)$ are computed using
the iteration \eqref{q:vin}--\eqref{q:chin}.
If this iteration
were convergent, say to $(\Vbb^\infty,\Chi^\infty,\Phi^\infty)$,
and if we let $\qa$ be the solution of
\begin{equation}
   \dy_t\qa + \Pl\bigl[ \sgb\cdot(\ub^\infty\cdot\gb\vb^\infty)
	- \dy_z(\ub^\infty\cdot\gb\rho^\infty) \bigr]
	= \mu \lapl3\qa + f_q^< \,,
\end{equation}
where $\ub^\infty$ and $\rho^\infty$ are determined from
$(\Vbb^\infty,\Chi^\infty,\Phi^\infty)$ as before,
then $(\vb^\infty,\rho^\infty)$ would be an exact solution of the OPEs.
As is often the case in this type of problems (see \cite{lorenz:80} and
\cite{jv-yavneh:04} for closely related problems), this iteration is asymptotic
rather than convergent, so we have to end the iteration at some $n=n_*$.
It is shown below that for $n=0,\cdots,n_*=\lfloor\eta/\eps^{1/4}\rfloor$:
(i)~the slaved variables $(\Vbb^n,\Chi^n,\Phi^n)$ are bounded
[see \eqref{q:fxhypo} below], and
(ii)~the QGE${}^n$ has an error of order $\eps^n$ locally in time
[see \eqref{q:Rest} below].

\medskip\noindent * We then show that with $\bar\vba^n=\Vbb^n(\qa;\eps)$,
$\xa^n=\Chi^n(\qa;\eps)$ and $\fa^n=\Phi^n(\qa;\eps)$ the solution of
\eqref{q:qgen} is bounded for $0\le t\le T_0$ with $T_0$ independent
of $\eps$.
(One can in fact show that the solution is bounded for all time, but
we shall not do so here.)

\medskip\noindent * Finally, we obtain bounds for $\hat W$ using the usual
Gronwall-type argument.
This proves the main theorem since $W^>$ is exponentially small
by Lemma~\ref{t:gev}.

\medskip
Following a common practice, we write $c$ for a generic constant which may
not be the same each time it appears;
the more important constants are numbered: $\cnst1$, $\cnst2$, etc.
Unless otherwise indicated,
these constants (assumed positive) may depend on $s$, $\sigma$ and
domain size ($L_i$), but not on $\eps$, $\kappa$, $n$, $\mu$ or the
initial data.

\subsection{Construction of QGE${}^n$}

Let $\Vbb^0=0$, $\Chi^0=0$ and $\Phi^0=0$.
From \eqref{q:vbn}--\eqref{q:rn} we then have
\begin{equation}\label{q:vrz}
   \vb^0 = \sgb\ilapl3\qt
   \qquad\hbox{and}\qquad
   \rho^0 = -\dy_z\ilapl3\qt.
\end{equation}
Putting these in \eqref{q:qgen} and using the fact that
$\sgb\cdot(\vb^0\cdot\gb\vb^0)-\dy_z(\vb^0\cdot\gb\rho^0)=\vb^0\cdot\gb\qt$,
we find
\begin{equation}
   \dy_t\qt + \Pl\bigl[ \vb^0\cdot\gb\qt ] = \mu\lapl3\qt + f_q^< \,,
\end{equation}
which is the quasi-geostrophic equation \eqref{q:qge}--\eqref{q:qg-inv}
for the truncated variable $\qt$.
This is our zeroth-order (truncated) quasi-geostrophic model, QGE${}^0$.

For $n\ge1$,
we define $\Vbb^n(\qa;\eps)$, $\Chi^n(\qa;\eps)$ and $\Phi^n(\qa;\eps)$
by the following iteration [cf.\ \eqref{q:pe-qxf}]
\begin{align}
  -&\frac1\eps (\Vbb_z^{n+1})^\perp =
	\Pl\bigl\{(\FD\Vbb_z^n)\Gyle(\qa,\Vbb^n,\Chi^n,\Phi^n)
	+ \dy_z\bigl(\overline{w^n\vb^n_z\vphantom{f}}\bigr)\bigr\}
	- \mu\dy_{zz}\Vbb_z^n - \dy_z\bar\fv, \label{q:vin}\\
   &\frac1\eps \lapl3\Phi^{n+1} = \Pl\bigl\{(\FD\lapl3\Chi^n)\,
		\Gyle(\qa,\Vbb^n,\Chi^n,\Phi^n) \notag\\
	&\qquad+ \lapl3\ilapl2 \gb\cdot(\uba^n\cdot\gb\vba^n)\bigr\}
	- \mu\Delta_3^2\Chi^n
	- \lapl3 f_\chi^\sle, \label{q:phin}\raisetag{3em}\\
  -&\frac1\eps \lapl3\Chi^{n+1} = \Pl\bigl\{(\FD\lapl3\Phi^n)\,
		\Gyle(\qa,\Vbb^n\Chi^n,\Phi^n) \notag\\
	&\qquad+ \dy_{zz}\ilapl2 \sgb\cdot(\uba^n\cdot\gb\vba^n)
	+\Pz\dy_z(\uba^n\cdot\gb\rba^n)\bigr\} - \mu\lapl3\Phi_{zz}^n
	- \dy_{zz} f_\phi^\sle, \label{q:chin}\raisetag{3em}
\end{align}
where $\FD$ denotes derivative (differential) with respect to $\qa$ and
\begin{align}\label{q:Gydef}
   \Gyle(\qt,\Vbb^n,\Chi^n,\Phi^n) &= \Pl \bigl[ \mu\lapl3 \qt + f_q
	- \sgb\cdot(\ub^n\cdot\gb\vb^n)
	+ \dy_z(\ub^n\cdot\gb\rho^n) \bigr]\\
	&=: \Gylh(\qt,\Vbb^n,\Chi^n,\Phi^n) + \mu\lapl3\qt + f_q^<\,.\notag
\end{align}
We recall that $(\ub^n,\rho^n)$ are given in terms of $(\Vbb^n,\Chi^n,\Phi^n)$
by \eqref{q:vbn}--\eqref{q:wn}.
This construction was presented, in general and formally, in \cite{wbsv:95},
which also proposed the use of a series expansion.

\medskip
Let $s>3/2$ be fixed (in the theorem we use $s=2$, but we keep the
general $s$ in the proof), and consider $\qa$ fixed for now.
We start by estimating the nonlinear terms in
$\Gyle(\qa,\Vbb^n,\Chi^n,\Phi^n)=:\Gy^n$ using \eqref{q:ur-qxf} and
\eqref{q:vbn}--\eqref{q:wn}:
For $\dy_z(\ub^n\cdot\gb\rho^n)$, we compute
\begin{equation}\begin{aligned}
   |\vb^n|_s &\le c\,(|\Vbb^n|_s + |\gb_2\psi^n|_s
		+ |\gb_2\chi^n|_s)\\
   &\le c\,(|\Vbb^n|_s + |\qt|_{s-1} + |\gb_3\Phi^n|_s + |\gb_3\Chi^n|_s),
\end{aligned}\end{equation}
\begin{equation}\begin{aligned}
   |\vb^n_z|_s &\le c\,(|\Vbb^n_z|_s + |\gb_2\psi^n_z|_s
		+ |\gb_2\chi^n_z|_s)\\
	&\le c\,(|\Vbb^n_z|_s + |\qt|_s + |\lapl3\Phi^n|_s + |\lapl3\Chi^n|_s),
\end{aligned}\end{equation}
\begin{equation}
   |w^n|_s \le c\,|w^n_z|_s \le c\,|\lapl3\Chi^n|_s\,.
\end{equation}
From these we have
\begin{equation}\begin{aligned}
   c\,|\ub^n|_s &\le |\Vbb^n|_s + |\qt|_{s-1} + |\lapl3\Chi^n|_s + |\gb_3\Phi^n|_s\,\\
   c\,|\ub^n_z|_s &\le |\Vbb^n_z|_s + |\qt|_s + |\lapl3\Chi^n|_s + |\lapl3\Phi^n|_s\,.
\end{aligned}\end{equation}
Putting these together with
\begin{equation}\label{q:bdrn}
   c\,|\gb_3\rho^n|_s \le |\qt|_s + |\lapl3\Phi^n|_s
    \qquad\textrm{and}\qquad
   c\,|\gb_3\rho^n_z|_s \le |\qt|_{s+1} + |\lapl3\Phi^n|_{s+1}
\end{equation}
gives us
\begin{equation*}\begin{aligned}
   |\ub^n_z\cdot\gb\rho^n|_s &\le c\,
       \bigl(|\qt|_s + |\Vbb^n_z|_s + |\lapl3\Phi^n|_s + |\lapl3\Chi^n|_s\bigr)
       \bigl(|\qt|_s + |\lapl3\Phi^n|_s\bigr)\\
   |\ub^n\cdot\gb\rho^n_z|_s &\le c\,
       \bigl(|\qt|_{s-1} + |\Vbb^n|_s + |\gb_3\Phi^n|_s +|\lapl3\Chi^n|_s\bigr)
       \bigl(|\qt|_{s+1} + |\lapl3\Phi^n|_{s+1}\bigr)\,,\\
\end{aligned}\end{equation*}
and, upon the use of \eqref{q:ipoi},
\begin{equation}\label{q:estrz1}
   |\dy_z(\ub^n_z\cdot\gb\rho^n)|_s
	\le c\kappa\,\bigl(|\qt|_s + |\Vbb^n_z|_s + |\lapl3\Phi^n|_s
		+ |\lapl3\Chi^n|_s\bigr)^2.
\end{equation}
For the other nonlinear term in $\Gy^n$,
\begin{equation}
   \sgb\cdot(\ub^n\cdot\gb\vb^n) = \ub^n\cdot\gb\lapl2\psi^n
	+ (\sgb\ub^n):\gb_3(\Vbb^n+\sgb\psi^n+\gb\Chi^n),
\end{equation}
we need in addition
\begin{equation}\begin{aligned}
   |\gb_2\ub^n|_s &\le c\,(|\lapl2\psi^n|_s + |\lapl2\Chi^n|_s + |\gb_2 w^n|_s)\\
	&\le c\,(|\qt|_s + |\lapl3\Chi^n|_{s+1} + |\lapl3\Phi^n|_s),
\end{aligned}\end{equation}
\begin{equation}\begin{aligned}\label{q:bdvn}
   |\gb_3\vb^n|_s &\le c\,(|\Vbb^n_z|_s + |\lapl3\psi^n|_s + |\lapl3\Chi^n|_s)\\
	&\le c\,(|\qt|_s + |\Vbb^n_z|_s + |\lapl3\Chi^n|_s + |\lapl3\Phi^n|_s),
\end{aligned}\end{equation}
and
\begin{equation}
   c\,|\gb_3\lapl2\psi^n|_s \le |\qt|_{s+1} + |\lapl3\Phi^n|_{s+1}\,.
\end{equation}
Putting these together and using \eqref{q:ipoi} again, we find
\begin{equation}\label{q:bdnlnv}
   |\sgb\cdot(\uba^n\cdot\gb\vba^n)|_s
	\le c\kappa\, \bigl(|\qa|_s + |\Vbb_z^n|_s + |\lapl3\Phi^n|_s
		+ |\lapl3\Chi^n|_s\bigr)^2.
\end{equation}
The linear terms are easily bounded, giving us
\begin{equation}
   |\Gy^n|_s \le c_1\kappa\, \bigl(|\qa|_s + |\Vbb_z^n|_s + |\lapl3\Phi^n|_s
		+ |\lapl3\Chi^n|_s\bigr)^2
		+ \mu\kappa^2 |\qa|_s + |f_q^\sle|_s.
\end{equation}

Going next to the nonlinear terms in \eqref{q:vin}--\eqref{q:chin},
we compute
\begin{equation}
   |\overline{w^n\vb_z^n}|_s \le c\,|w^n|_s|\sgb\psi^n_z+\gb\Chi^n_z|_s
	\le c\kappa\, |\lapl3\Chi^n|_s
	\bigl( |\qa|_s + |\lapl3\Phi^n|_s + |\lapl3\Chi^n|_s \bigr).
\end{equation}
The following estimates follow from our computation above:
\begin{equation}\label{q:estvz1}\begin{aligned}
   |\lapl3\ilapl2\gb\cdot(\ub^n\cdot\gb\vb^n)|_s
	&\le c\kappa^3\, \bigl(|\qa|_s + |\Vbb_z^n|_s + |\lapl3\Phi^n|_s
		+ |\lapl3\Chi^n|_s\bigr)^2,\\
   |\dy_{zz}\ilapl2\sgb\cdot(\uba^n\cdot\gb\vba^n)|_s
	&\le c\kappa^3\, \bigl(|\qa|_s + |\Vbb_z^n|_s + |\lapl3\Phi^n|_s
		+ |\lapl3\Chi^n|_s\bigr)^2.
\end{aligned}\end{equation}

Using the relations \eqref{q:fuvr-fqxf}, the forcing terms are bounded as
\begin{equation}\label{q:fqxf}\begin{aligned}
   &|f_q|_s^{} \le c \bigl(|\fv|_{s+1}^{} + |\fr|_{s+1}^{}\bigr)
   &\qquad &|\lapl3 f_\chi|_s^{} \le c |\fv|_{s+3}^{}\\
   &|\dy_z\bar\fv|_s^{} \le |\fv|_{s+1}^{}
   &&|\dy_{zz}f_\phi|_s^{} \le c \bigl(|\fv|_{s+3}^{} + |\fr|_{s+3}^{}\bigr).
\end{aligned}\end{equation}

\medskip
Moving on to terms such as $\FD(\lapl3\Chi^n)\Gy^n$ in
\eqref{q:vin}--\eqref{q:chin}, we need a few definitions:
Let $D_\eta(\qa)$ be the complex $\eta$-neighbour\-hood of $\qa$
in $\Pl H^s(\Dom)$.
It can be defined by Fourier series as
\begin{equation}\label{q:deta-def}
   D_\eta(\qa) = \Bigl\{ q' : q'(\xb) = \sum_{|\kb|<\kappa}
	q_\kb' \e^{{\rm i}\kb\cdot\xb} \textrm{ with }
	\sum_{|\kb|<\kappa} |\kb|^{2s}\,|q_\kb'-\qa_\kb^{}|^2<\eta^2 \Bigr\}.
\end{equation}
We note that since $\qt(\xb)$ is real, the Fourier coefficients must
satisfy $\qt_{-\kb}=\overline{\qt_{\kb}}$ where (here and only here) bar
denotes complex conjugate;
however, the Fourier coefficients $q_\kb'$ in \eqref{q:deta-def} need not
satisfy this constraint.
Let $\delta>0$ be fixed.
For any function $g$ of $\qa$, let
\begin{equation}
   |g|_{s;n}^{} := \sup_{D_{\eta-n\delta}(\qa)} |g(\qa)|_s^{},
\end{equation}
which is meaningful for $n=0,\cdots,\lfloor\eta/\delta\rfloor=:n_*$
when $D_{\eta-n\delta}(\qa)$ is non-empty.

Let us also denote
\begin{equation}
   |U^n|_{s;n} := |\Vbb_z^n|_{s;n} + |\lapl3\Phi^n|_{s;n}
	+ |\lapl3\Chi^n|_{s;n}.
\end{equation}
To bound  $\FD(\Vbb_z^n)\Gy^n$, $\FD(\lapl3\Chi^n)\Gy^n$ and
$\FD(\lapl3\Phi^n)\Gy^n$, we use Cauchy's integral formula
(cf.\ e.g., \cite{lochak-meunier:mpavg}):
For $\varphi:{\rm cl}\,D_\eta(\qa)\to\Comp$ analytic with
${\rm cl}\,D_\eta(\qa)\subset\Comp^l$ and
$\delta>0$, one can bound $|\varphi'|$ in ${\rm cl}\,D_{\eta-\delta}(\qa)$ by
$|\varphi|$ in ${\rm cl}\,D_\eta(\qa)$,
\begin{equation}
   |\varphi'\cdot z|_{{\rm cl}\,D_{\eta-\delta}(\qa)} \le
	\frac1\delta |\varphi|_{{\rm cl}\,D_\eta(\qa)} |z|_{\Comp^l}.
\end{equation}
Now $\Vbb_z(\qa;\eps)$, $\lapl3\Chi^n(\qa;\eps)$ and $\lapl3\Phi^n(\qa;\eps)$
are analytic functions of $\qa$, being polynomials in $\qa_\kb$.
Using this, we have
\begin{equation}\label{q:fdfn}\begin{aligned}
   |(\FD\lapl3\Phi^n)\Gy^n|_{s,n+1}^{}
	&\le \frac{c}\delta|\lapl3\Phi^n|_{s,n}^{}|\Gy^n|_{s,n}^{}\\
	&\le \frac{c}\delta|\lapl3\Phi^n|_{s,n}^{}
		\bigl\{ \mu\kappa^2|\qa|_s + |f_q^<|_s
		+ c\kappa^2\,\bigl( |\qa|_s + |U^n|_{s;n} \bigr)^2 \bigr\}.
\end{aligned}\end{equation}

We have now bounded every term in \eqref{q:chin}: putting together
\eqref{q:fdfn}, \eqref{q:fqxf}, \eqref{q:estvz1} and \eqref{q:estrz1},
we find
\begin{equation}\begin{aligned}
   \frac1\eps |\lapl3\Chi^{n+1}|_{s;n+1} &\le
	\mu\kappa^2\,|\lapl3\Phi^n|_{s;n} + |\dy_{zz}f_\phi^<|_s
	+ c\kappa^3\, \bigl( |\qa|_s + |U^n|_{s;n} \bigr)^2\\
	&\quad+ \frac{c}{\delta}|\lapl3\Phi^n|_{s;n}
	    \bigl\{ \mu\kappa^2|\qa|_s + |f_q^<|_s
		+ c\kappa^2\,\bigl( |\qa|_s + |U^n|_{s;n} \bigr)^2 \bigr\}.
\end{aligned}\end{equation}
Similarly, we find for \eqref{q:phin},
\begin{equation}\begin{aligned}
   \frac1\eps |\lapl3\Phi^{n+1}|_{s;n+1} &\le
	\mu\kappa^2\,|\lapl3\Chi^n|_{s;n} + |\lapl3 f_\chi^<|_s
	+ c\kappa^3\, \bigl( |\qa|_s + |U^n|_{s;n} \bigr)^2\\
	&\quad+ \frac{c}{\delta}|\lapl3\Phi^n|_{s;n}
	    \bigl\{ \mu\kappa^2|\qa|_s + |f_q^<|_s
		+ c\kappa^2\,\bigl( |\qa|_s + |U^n|_{s;n} \bigr)^2 \bigr\},
\end{aligned}\end{equation}
and for \eqref{q:vin},
\begin{equation}\begin{aligned}
   \frac1\eps |\Vbb^{n+1}_z|_{s;n+1} &\le
	\mu\kappa^2\,|\Vbb_z^n|_{s;n} + |\dy_{z}\bar\fv^<|_s
	+ c\kappa^3\, \bigl( |\qa|_s + |U^n|_{s;n} \bigr)^2\\
	&\quad+ \frac{c}{\delta}|\lapl3\Phi^n|_{s;n}
	    \bigl\{ \mu\kappa^2|\qa|_s + |f_q^<|_s
		+ c\kappa^2\,\bigl( |\qa|_s + |U^n|_{s;n} \bigr)^2 \bigr\}.
\end{aligned}\end{equation}
Putting these estimates together with \eqref{q:fqxf}, we find
\begin{equation}\label{q:uiter}\begin{aligned}
   \frac1\eps |U^{n+1}|_{s;n+1} &\le
	\mu\kappa^2\,|U^n|_{s;n} + \|f\|
	+ \cnst2\kappa^3\, \bigl( |\qa|_s + |U^n|_{s;n} \bigr)^2\\
	&\qquad+ \frac{\cnst3}{\delta}\,|U^n|_{s;n}
	    \bigl\{ \mu\kappa^2|\qa|_s + |f_q^<|_s
	    + \cnst4\kappa^2\,\bigl( |\qa|_s + |U^n|_{s;n} \bigr)^2 \bigr\},
\end{aligned}\end{equation}
where $\|f\|:=|(f_{\vb}^{<},f_{\rho}^{<})|_{s+3}^{}$.

\medskip
Now suppose that at iteration $n$, with $\cnst1=8\cnst2$ we have
\begin{equation}\label{q:fxhypo}
   |U^n|_{s;n}^{} \le \cnst{1}\, \eps^{1/4}\,
	\bigl( |\qa|_s^2 + \|f\| \bigr).
\end{equation}
Since $\Vbb^0=\Chi^0=\Phi^0=0$, this trivally holds for $n=0$.
Taking $\delta=\eps^{1/4}$ and $\kappa=\eps^{-1/4}$,
and substituting this into \eqref{q:uiter}, we find that
\begin{equation}
   |U^{n+1}|_{s;n+1}^{} \le \cnst{1}\, \eps^{1/4}\,
	\bigl( |\qa|_s^2 + \|f\| \bigr).
\end{equation}
provided that $\eps$ satisfies the following constraints
\begin{equation}\label{q:eps1}\begin{aligned}
   \eps \le \min\Bigl(\, 1/(16\mu^2),\, &1/c_1^4,\, (\cnst1/4)^{4/3},\,
	1/(12\mu\cnst3|\qa|_s^{})^4,\,\\ &1/(12\cnst3\|f\|)^{4/3},\,
	1/\bigl(12\cnst3\cnst4[2|\qa|_s^{}+\|f\|]^2\bigr)^4 \,\Bigr).
\end{aligned}\end{equation}
By induction \eqref{q:fxhypo} can be seen to hold for
$n=1,\cdots,n_*$.
For future reference, we note that this restriction on $\eps$ and
\eqref{q:fxhypo} imply that
\begin{equation}\label{q:bdqU}
  |\qa|_s^{} + |U^n|_s^{} \le c_2^{} (|\qa|_s^{} + \|f\|) \,.
\end{equation}

\subsection{Local Error Bounds}

Moving on to estimating the approximation error locally in time, we define
\begin{align}
   \Rv^n &:= (\FD\Vbb_z^n)\Gy^n + \frac1\eps \Vbb_z^n
	+ \dy_z(\overline{w^n\vb_z^n})
	- \mu\dy_{zz}^{}\Vbb_z^n - \dy_z \bar{f}^<_{\vb} \label{q:Rvdef}\\
	&= \frac1\eps \bigl[ \Vbb_z^n - \Vbb_z^{n+1} \bigr]\notag\\
   \Rx^n &:= (\FD\lapl3\Chi^n)\Gy^n
	- \frac1\eps \lapl3\Phi^n + \lapl3\ilapl2\gb\cdot(\ub^n\cdot\gb\vb^n)
	- \mu\Delta_3^2\Chi^n - \lapl3 f_\chi^< \\
	&= \frac1\eps \bigl[ \lapl3\Phi^{n+1} - \lapl3\Phi^n \bigr]  \notag\\
   \Rf^n &:= (\FD\Phi_{zz}^n)\Gy^n
	+ \frac1\eps \lapl3\Chi^n
	+ \dy_{zz}\ilapl2\sgb\cdot(\ub^n\cdot\gb\vb^n) \label{q:Rfdef}\\
	&\qquad+ \Pz\dy_z(\ub^n\cdot\gb\rho^n)
	- \mu\lapl3\Phi_{zz}^n - \dy_{zz} f_\phi^<   \notag\\
	&= \frac1\eps \bigl[ \lapl3\Chi^n - \lapl3\Chi^{n+1} \bigr]. \notag
\end{align}

Since $\Vbb^0=0$, $\Chi^0=0$ and $\Phi^0=0$, we have
\begin{align}
   \Rv^0 &= -\dy_z\bar{f}^<_{\vb}, \notag\\
   \Rx^0 &= \lapl3\ilapl2\gb\cdot(\vb^0\cdot\gb\vb^0) 
	- \lapl3 f_\chi^<, \notag\\
   \Rf^0 &= \dy_{zz}\ilapl2\sgb\cdot(\vb^0\cdot\gb\vb^0)
	+ \Pz\dy_z(\vb^0\cdot\gb\rho^0) - \dy_{zz} f_\phi^<, \notag
\end{align}
which leads to the estimate, denoting
$|\Rr^n|_s^{} := |\Rv^n|_s^{} +|\Rx^n|_s^{} +|\Rf^n|_s^{}$,
\begin{equation}
   |\Rr^0|_s^{} \le C\eps^{-3/4}\,\bigl(|\qt|_s^{} + \|f\|\bigr)^2 \,.
\end{equation}

For $n\ge1$,
let us denote $\dl\vb^n:=\vb^{n+1}-\vb^n$, $\dl\rho^n:=\rho^{n+1}-\rho^n$,
etc., and consider
\begin{equation}\label{q:rvn}\begin{aligned}
   \Rv^{n+1} &= (\FD\Vbb_z^{n+1})\Gy^{n+1} + \frac1\eps\Vbb_z^{n+1}
	+ \dy_z(\overline{w^{n+1}\vb_z^{n+1}})
	- \mu\dy_{zz}\Vbb_z^{n+1} - \dy_z \bar{f}^<_{\vb}\\
	&= (\FD\Vbb_z^{n+1})\dl\Gy^n - \eps (\FD\Rv^n)\Gy^n\\
	&\qquad+ \dy_z(\overline{\dl w^n\vbb_z^n})
	+ \dy_z(\overline{w^{n+1}\dl\vb_z^n}) + \eps\mu\dy_{zz}\Rv^n\,,
\end{aligned}\end{equation}
For the second equality, we have used the second equalities in
\eqref{q:Rvdef}--\eqref{q:Rfdef} to write $\dl\vbb^n$ in terms of $\Rv^n$,
etc.
To estimate these $\dl$-quantities, we compute
\begin{equation}\begin{aligned}
    &\lapl3\delta\rho^n = \lapl2[\dy_z\delta\Phi^n]\\
    \Rightarrow\> |&\gb_3\dl\rho^n|_s^{} \le c\,|\lapl3\delta\Phi^n|_s^{}
    \le \eps c\,|\Rx^n|_s^{}\,.
\end{aligned}\end{equation}
Similarly, we find
\begin{equation}\begin{aligned}
   &|\dy_z\dl w^n|_s^{} = |\lapl2\dl\Chi^n|_s
	\le \eps c\, |\Rf^n|_s^{}\\
   &|\gb_3\dl\vb^n|_s^{} \le c\,|\gb_3\dl\Vbb^n|_s
	+ |\lapl3\dl\psi^n|_s + |\lapl3\dl\Chi^n|_s\\
   &\phantom{|\gb_3\dl\vb^n|_s^{}} \le c\,\bigl( |\dl\Vbb^n_z|_s
		+ |\dy_{zz}\dl\Phi^n|_s + |\lapl3\dl\Chi^n|_s \bigr)
	\le \eps c\, |\Rr^n|_s^{} \\
   \Rightarrow\quad &|\dy_z\dl\ub^n|_s^{} \le \eps c\, |\Rr^n|_s^{} \,.
\end{aligned}\end{equation}
Estimating the terms in
\begin{equation}\begin{aligned}
   \dl\Gy^n &= \dy_z(\dl\ub^n\cdot\gb\rho^n + \ub^{n+1}\cdot\gb\dl\rho^n)
	- \sgb\cdot(\dl\ub^n\cdot\gb\vb^n + \ub^{n+1}\cdot\gb\dl\vb^n)\\
	&= \dl\ub^n_z\cdot\gb\rho^n + \dl\ub^n\cdot\gb\rho^n_z
	+ \ub^{n+1}_z\cdot\gb\dl\rho^n + \ub^{n+1}\cdot\gb\dl\rho^n_z\\
	&\quad - \sgb\dl\ub^n:\gb\vb^n - \dl\ub^n\cdot\gb\lapl2\psi^n
	- \sgb\ub^{n+1}:\gb\dl\vb^n - \ub^{n+1}\cdot\gb\lapl2\dl\psi^n\\
\end{aligned}\end{equation}
using what we have computed above as
\begin{equation*}\begin{aligned}
   c\,|\dl\Gy^n|_s &\le \eps\,|\Rr^n|_s(|\qt|_s+|U^n|_s)
	+ \eps\kappa\,|\Rr^n|_s(|\qt|_s+|U^n|_s)
	+ \eps\,(|\qt|_s+|U^n|_s)|\Rr^n|_s\\
	&\>+ \eps\kappa\,(|\qt|_s+|U^n|_s)|\Rr^n|_s
	+ \eps\kappa\,|\Rr^n|_s(|\qt|_s+|U^n|_s)
	+ \eps\,|\Rr^n|_s(|\qt|_s+|U^n|_s)\\
	&\>+ \eps\kappa\,(|\qt|_s+|U^n|_s)|\Rr^n|_s
	+ \eps\kappa\,(|\qt|_s+|U^n|_s)|\Rr^n|_s\,,
\end{aligned}\end{equation*}
we obtain
\begin{equation}
   |\dl\Gy^n|_s^{} \le \eps^{3/4} c\,|\Rr^n|_s^{}
	\bigl( |\qa|_s^{} + |U^n|_s^{} \bigr)
	\le \eps^{3/4} c\,|\Rr^n|_s^{} \bigl( |\qa|_s^{} + \|f\| \bigr).
\end{equation}

Estimating $(\FD\Vbb^{n+1}_z)\dl\Gy^n$ and $(\FD\Rv^n)\Gy^n$ in \eqref{q:rvn}
using Cauchy's formula as before and using \eqref{q:bdqU}, we find
\begin{equation}\begin{aligned}
   c\,|\Rv^{n+1}|_{s;n+1}^{} &\le
	\frac{1}{\dl}\,|\Vbb_z^n|_{s;n}^{}|\dl\Gy^n|_s^{}
	+ \frac{\eps}{\dl}\,|\Rv^n|_{s;n}|\Gy^n|_s\\
     &\qquad+ |\dy_z(\overline{\dl w^n\vbb_z^n})|_s^{}
	+ |\dy_z(\overline{w^{n+1}\dl\vb_z^n})|_s^{}
	+ \eps\mu\kappa^2|\Rv^n|_{s;n}^{}\\
     &\le \frac{\eps\kappa}{\dl}\,|\Vbb^n_z|_{s;n}|\Rr^n|_s(|\qa|_s+\|f\|)
	+ \frac{\eps\kappa}{\dl}\,(|\qa|_s^{} + \|f\|)^2|\Rv^n|_{s;n}^{}\\
	&\qquad+ \frac{\mu\eps\kappa^2}{\dl}\,|\Rv^n|_{s;n}^{}\,|\qa|_s^{}
	+ \eps|\Rr^n|(|\qa|_s + \|f\|) + \mu\eps\kappa^2 |\Rr^n|_s\\
     |\Rv^{n+1}|_{s;n+1}^{} &\le c_5'\, \eps^{1/4}
	\bigl( |\qa|_s^{} + \|f\| + \cnst6(\mu) \bigr)^2
	|\Rr^n|_{s;n}
\end{aligned}\end{equation}
The estimates for $\Rx^n$ and $\Rf^n$ are identical in form, giving us
\begin{equation}
   |\Rr^{n+1}|_{s;n+1}^{} \le \cnst5\, \eps^{1/4}
	\bigl( |\qa|_s^{} + \|f\| + \cnst6 \bigr)^2
	|\Rr^n|_{s;n} \,.
\end{equation}
So $|\Rr^n|_{s;n}^{}$ is decreasing in $n$ for $\eps$ sufficiently small.
Remembering that this construction is valid for $n=0,\cdots,n_*$,
$|\Rr^n|_{s;n}^{}$ is therefore smallest for $n=n_*=\eta\eps^{-1/4}$.
More precisely, when
\begin{equation}\label{q:eps2}
    \cnst5\, \eps^{1/4} \bigl(|\qa|_s^{} + \|f\| + \cnst6\bigr)^2
	\le K < 1\,,
\end{equation}
we have
\begin{equation}\label{q:Rest0}\begin{aligned}
   |\Rr^{n_*}|_{s;n_*}^{} &\le  |\Rr^0|_{s;0}^{} \, \Bigl[ \cnst5 \eps^{1/4}
		\bigl(|\qa|_s^{} + \|f\| + \cnst6\bigr)^2 \Bigr]^{n_*}\\
	&\le |\Rr^0|_{s;0}^{}
	   \exp\bigl( \eps^{-1/4}\,\eta\log K \bigr)\\
	&\le c\, \bigl(|\qa|_s^{} + \|f\|\bigr)^2
	\eps^{-3/4}\,\exp\bigl( \eps^{-1/4}\,\eta\log K \bigr)
\end{aligned}\end{equation}
which is exponentially small in $\eps$ (the argument of the exponential
is negative since $K<1$).
In what follows we fix $K$, say $K=1/2$,
and set $\eta=\sigma/\log 2$.
Assuming in addition that $\eps\le\sigma/10$, we can write
$\eps^{-3/4}\exp(-\sigma/\eps^{1/4}) \le c\,\exp(-\sigma/\eps^{1/4})$
and
\begin{equation}\label{q:Rest}
   |\Rr^{n_*}|_{s;n_*}^{} \le c\, \bigl(|\qa|_s^{} + \|f\|\bigr)^2
	\,\exp( -\sigma/\eps^{-1/4} ).
\end{equation}

\medskip\noindent{\bf Remarks.}
A similar result for a different model was obtained in \cite{dw:abal}.
This result, which does not depend on Gevrey regularity (the model
was inviscid), appears to be peculiar to the model in question and is
only local in time (i.e.\ only up to this point in the present proof).


\subsection{Regularity of QGE${}^n$}

We now show that the solution $\qa(t)$ of
\begin{equation}\label{q:dqadt}
   \dy_t\qa + \Pl\bigl[ \sgb\cdot(\ub^*\cdot\gb\vb^*)
	- \dy_z(\ub^*\cdot\gb\rho^*) \bigr] = \mu \lapl3\qa + f_q^<
\end{equation}
where $\ub^*=\ub^{n_*}$ and $\rho^*=\rho^{n_*}$, is bounded independently
of $\eps$ for $t\in[0,T_0]$.
To this end we write $\ub^*=\ub^0+\ub^\eps$ and $\rho^*=\rho^0+\rho^\eps$
where $\ub^0=(u^0,v^0,0)$,
\begin{equation}
   \vb^0 = \sgb\ilapl3\qa
   \qquad\textrm{and}\qquad
   \rho^0 = -\dy_z\ilapl3\qa.
\end{equation}
As a preparation, we note that (\ref{q:bdrn}a) and \eqref{q:bdvn} imply,
with \eqref{q:fxhypo},
\begin{equation}\begin{aligned}
   &|\gb_3\rho^\eps|_s \le c\,|\lapl3\Phi^*|_s
	\le c\,\eps^{1/4}\,(|\qa|_s^2 + \|f\|)\\
   &|\gb_3\vb^\eps|_s \le c\,\bigl(|\Vbb^*_z|_s + |\lapl3\Chi^*|_s
	+ |\lapl3\Phi^*|_s\bigr) \le c\,\eps^{1/4}\,(|\qa|_s^2 + \|f\|).\\
\end{aligned}\end{equation}
Noting that
$\sgb\cdot(\vb^0\cdot\gb\vb^0)-\dy_z(\vb^0\cdot\gb\rho^0)=\vb^0\cdot\gb\qa$,
we write \eqref{q:dqadt} as
\begin{equation}\begin{aligned}
   &\dy_t\qa + \Pl\bigl[ \vb^0.\gb\qa + \sgb\cdot(\ub^\eps\cdot\gb\vb^0)
     + \sgb\cdot(\vb^0\cdot\gb\vb^\eps) + \sgb\cdot(\ub^\eps\cdot\gb\vb^\eps)\\
     &\qquad- \dy_z(\ub^\eps\cdot\gb\rho^0) - \dy_z(\vb^0\cdot\gb\rho^\eps)
       - \dy_z(\ub^\eps\cdot\gb\rho^\eps) \bigr] = \mu \lapl3\qa + f_q^<.
\end{aligned}\end{equation}

We multiply this by $\Delta_3^s\qa$ in $L^2(\Dom)$ and estimate
the resulting terms as
\begin{equation}\begin{aligned}
   |(\Delta_3^s\qa,\vb^0\cdot\gb\qa)|
	&\le c\,|\qa|_s^{}|\vb^0\cdot\gb\qa|_s^{}
	\le c\,|\qa|_s^{} |\qa|_{s-1}^{} |\qa|_{s+1}^{}\\
	&\le \frac\mu4 |\qa|_{s+1}^2 + \frac{c}\mu |\qa|_{s-1}^2 |\qa|_s^2\,;
\end{aligned}\end{equation}
and for the terms involving $\vb$,
\begin{equation}\begin{aligned}
   |(\Delta_3^s\qa,\sgb\cdot(\ub^\eps\cdot\gb\vb^0))|
	&\le c\,|\qa|_s|\sgb\cdot(\ub^\eps\cdot\gb\vb^0)|_s
	\le c\,|\qa|_s\,\kappa|\ub^\eps\cdot\gb\vb^0|_s\\
	&\le c\,|\qa|_s\,\kappa|\ub^\eps|_s|\gb_3\vb^0|_s
	\le c\,|\qa|_s\,\kappa\eps^{1/4}(|\qa|_s^2+\|f\|)\,|\qa|_s\\
	&\le c\,|\qa|_s^2\,(|\qa|_s^2+\|f\|),\\
\end{aligned}\end{equation}
\begin{equation}
   |(\Delta_3^s\qa,\sgb\cdot(\vb^0\cdot\gb\vb^\eps))|
	\le c\,|\qa|_s^2\,(|\qa|_s^2+\|f\|),
\end{equation}
\begin{equation}\begin{aligned}
   |(\Delta_3^s\qa,\sgb\cdot(\ub^\eps\cdot\gb\vb^\eps))|
	&\le c\,|\qa|_s\,\kappa|\ub^\eps|_s\,|\gb_3\vb^\eps|_s\\
	&\le c\,\eps^{1/4}\,|\qa|_s\,(|\qa|_s^2+\|f\|);\\
\end{aligned}\end{equation}
and for those terms involving $\rho$,
\begin{equation}
   |(\Delta_3^s\qa,\dy_z(\ub^\eps\cdot\gb\rho^0)|
	\le c\,|\qa|_s^2\,(|\qa|_s^2+\|f\|),
\end{equation}
\begin{equation}
   |(\Delta_3^s\qa,\dy_z(\vb^0\cdot\gb\rho^\eps)|
	\le c\,|\qa|_s^2\,(|\qa|_s^2+\|f\|),
\end{equation}
\begin{equation}
   |(\Delta_3^s\qa,\dy_z(\ub^\eps\cdot\gb\rho^\eps)|
	\le c\,\eps^{1/4}\,|\qa|_s\,(|\qa|_s^2+\|f\|)^2.
\end{equation}
Putting these together, we have
\begin{equation}\label{q:dqasdt}
   \ddt{}|\qa|_s^2 + \mu\,|\qa|_{s+1}^2 \le
	\cnst8(\mu)\, |\qa|_s^{} \bigl(|\qa|_s^2 + \|f\| + c'\bigr)^2
	+ c''(\mu) \|f\|^2,
\end{equation}
where $\cnst8$ and $c'$ are independent of $\eps$ if one assumes
\eqref{q:eps1} and \eqref{q:eps2}.
Now let $Q=2\,|\qa(0)|_s^{}+\|f\|$.
Replacing $|\qa|_s^{}$ by $Q$ in \eqref{q:eps1} and \eqref{q:eps2},
we assume
\begin{equation}\label{q:epsQ}\begin{aligned}
   \eps \le \eps_0 = \min\Bigl(\, &\sigma/10,\, 1/(16\mu^2),\, 1/c_1^4,\,(\cnst1/4)^{4/3},\,
	1/(12\mu\cnst3 Q)^4,\,\\ 
	&\quad 1/\bigl[2\cnst5 \bigl(Q+\|f\|+\cnst6\bigr)\bigr]^4,\,
	1/(12\cnst3\|f\|)^{4/3},\\
	&\quad 1/\bigl[12\cnst3\cnst4(2Q+\|f\|)^2\bigr]^4 \,\Bigr).
\end{aligned}\end{equation}
It now follows from \eqref{q:dqasdt} that there exists a
$T_0(|\qa(0)|_s^{},\|f\|;\mu,s,\sigma,\Dom)>0$ such that for $t\in[0,T_0]$,
\begin{equation}\label{q:qabd}
   |\qa(t)|_s^{} \le 2\,|\qa(0)|_s^{}+\|f\|.
\end{equation}

\noindent{\bf Remarks.}
It is clear from the foregoing that this result is also valid for
any $n\le n_*$.
When $n=0$ (i.e.\ the classical quasi-geostrophic equation), one does not
need the $\kappa$ cutoff to prove boundedness (cf.\ \cite{bourgeois-beale:94}),
but this seems unavoidable for $n>0$.


\subsection{Long-time Error Bounds}

Let $W^*=(\vb^*,\rho^*)=(\vb^{n_*},\rho^{n_*})$.
We now seek to show that $\hat W = W^< - W^*$, {\em if initially
exponentially small, remains exponentially small over a timescale
of order one\/}.
We use the usual Gronwall-type argument.
\vfill\eject

Starting with the evolution equation for $\rho^<$,
\begin{equation}\begin{aligned}
   \dy_t\rho^< - \frac1\eps w^< + \Pl(\ub^<\cdot\gb\rho^<)
	&= \mu\lapl3\rho^< + f_\rho^<\\
	&\quad - \Pl\bigl(\ub^<\cdot\gb\rho^> + \ub^>\cdot\gb\rho^<
		+ \ub^>\cdot\gb\rho^>\bigr)\\
	&=: \mu\lapl3\rho^< + f_\rho^< + \Fy_\rho^>,
\end{aligned}\end{equation}
we find that $\hat\rho=\rho^*-\rho^<$ is governed by
\begin{equation}\label{q:drhdt}
   \dy_t\hat\rho - \frac1\eps \hat w + \Pl\bigl(\hat\ub\cdot\gb\hat\rho
	+ \hat\ub\cdot\gb\rho^* + \ub^*\cdot\gb\hat\rho\bigr)
	= \Sy_\rho^* + \mu\lapl3\hat\rho + \Fy_\rho^>,
\end{equation}
where
\begin{equation}
	\Sy_\rho^* = -\dy_t\rho^* + \frac1\eps w^* - \Pl(\ub^*\cdot\gb\rho^*)
		+ \mu\lapl3\rho^* + f_\rho^<.
\end{equation}
Now $\Sy_\rho^*$ can be expressed in terms of $\Rf^*$ as follows:
\begin{equation}\begin{aligned}
   -\lapl3\dy_z\Sy_\rho^* &= \dy_t\bigl(\lapl2\Phi^*_{zz} - \dy_{zz}\qa\bigl)
	- \frac{1}{\eps}\lapl3\dy_z w^* + \lapl3\dy_z(\ub^*\cdot\gb\rho^*)\\
	&\qquad- \mu\lapl3\dy_z\rho^* - \lapl3\dy_z f_\rho^<\\
	&= \lapl2(\FD\Phi^*_{zz})\,\Gy^* - \dy_{zz}\Gy^*
	+ \frac1\eps\lapl3\lapl2\Chi^* + \lapl3\dy_z(\ub^*\cdot\gb\rho^*)\\
	&\qquad- \mu\lapl3\dy_z\rho^* - \lapl3\dy_z f_\rho^<\,,
\end{aligned}\end{equation}
which gives upon using the definitions \eqref{q:Gydef} of $\Gy^*$ and
\eqref{q:Rfdef} of $\Rf$
\begin{equation}\label{q:Syr}
   \lapl3\dy_z\Sy_\rho^* = \lapl2\Rf^*.
\end{equation}

Similarly, we find for $\hat\vb$,
\begin{equation}\label{q:dvhdt}\begin{aligned}
   \dy_t\hat\vb + \frac1\eps \bigl[ \hat\vb^\perp + \gb\hat p \bigr]
	+ \Pl\bigl( \hat\ub\cdot\gb\hat\vb + \hat\ub\cdot\gb\vb^*
		&+ \ub^*\cdot\gb\hat\vb \bigr)\\
   &= \mu\lapl3\hat\vb + \Sy_{\vb}^* + \Fy_{\vb}^>
\end{aligned}\end{equation}
where
\begin{equation}
   \Fy_{\vb}^> = -\Pl\bigl(\ub^<\cdot\gb\vb^> + \ub^>\cdot\gb\vb^<
		+ \ub^>\cdot\gb\vb^> \bigr)
\end{equation}
and, following the computation leading to \eqref{q:Syr},
\begin{equation}\label{q:Syv}
   \bar\Sy_{\vb}^* = \Rv^* \,,\qquad
   \lapl3\sgb\cdot\Sy_{\vb}^* = \lapl2\Rf^*  \qquad\hbox{and}\qquad
   \lapl3\gb\cdot\Sy_{\vb}^* = \lapl2\Rx^* \,.
\end{equation}

Multiplying \eqref{q:dvhdt} by $\hat\vb$ and \eqref{q:drhdt} by $\hat\rho$
in $L^2(\Dom)$, we find
\begin{equation}\label{q:dbdvhdt}\begin{aligned}
   \frac12&\ddt{} \bigl( |\hat\vb|^2 + |\hat\rho|^2 \bigr)
	+ \mu \bigl( |\gb_3^{}\hat\vb|^2 + |\gb_3^{}\hat\rho|^2 \bigr)\\
	&= -(\hat\vb,\hat\ub\cdot\gb\vb^*) - (\hat\rho,\hat\ub\cdot\gb\rho^*)
	+ (\hat\vb,\Sy_{\vb}^*) + (\hat\rho,\Sy_\rho^*)
	+ (\hat\vb,\Fy_{\vb}^>) + (\hat\rho,\Fy_\rho^>).
\end{aligned}\end{equation}
Noting that
\begin{equation}\begin{aligned}
   &|\gb_3\vb^*|_{L^\infty(\Dom)}^{} \le c\,|\gb_3\vb^*|_s^{}
	\le c\,\bigl(|\qa|_s^{} + |U^*|_s^{}\bigr)
	\le c\,(|\qa|_s^{} + \|f\|)\,\\
   &|\gb_3\rho^*|_{L^\infty(\Dom)}^{} \le c\,|\gb_3\rho^*|_s^{}
	\le c\,\bigl(|\qa|_s^{} + |U^*|_s^{}\bigr)
	\le c\,(|\qa|_s^{} + \|f\|)\,,\\
\end{aligned}\end{equation}
where we have used \eqref{q:bdqU} for the last inequalities,
we bound the first two terms on the rhs of \eqref{q:dbdvhdt} by
\begin{equation}\begin{aligned}
   &|(\hat\vb,\hat\ub\cdot\gb\vb^*)| \le c\,|\hat\vb|\,|\gb_3^{}\hat\vb|\,
		|\gb_3^{}\vb^*|_{L^\infty(\Dom)}^{}
	\le \frac\mu6|\gb_3^{}\hat\vb|^2
		+ \frac{c}\mu |\hat\vb|^2\,(|\qa|_s^2 + \|f\|^2)\,,\\
   &|(\hat\rho,\hat\ub\cdot\gb\rho^*)| \le c\,|\hat\rho|\,|\gb_3^{}\hat\vb|\,
		|\gb_3^{}\rho^*|_{L^\infty(\Dom)}^{}
	\le \frac\mu6 |\gb_3^{}\hat\vb|^2
		+ \frac{c}\mu |\hat\rho|^2\,(|\qa|_s^2 + \|f\|^2)\,,
\end{aligned}\end{equation}
and use the Cauchy--Schwarz inequality for the remaining terms.
This gives us
\begin{equation}\label{q:mudvhdt}
   \mu\ddt{} \bigl( |\hat\vb|^2 + |\hat\rho|^2 \bigr)
	\le c\,(|\qa|_s^2 + \|f\|^2)\,\bigl(|\hat\vb|^2 + |\hat\rho|^2\bigr)
	+ c'\, |\Sy^*|^2 + c''\, |\Fy^>|^2,
\end{equation}
where we have denoted $|\Sy^*|^2 = |\Sy_{\vb}^*|^2 + |\Sy_{\rho}^*|^2$
and $|\Fy^>|^2 = |\Fy_{\vb}^>|^2 + |\Fy_\rho^>|^2$.

By \eqref{q:qabd}, we can replace $|\qa(t)|_s^{}$ by
$2\,|\qa(0)|_s^{} + \|f\|$ for $t\in[0,T_0]$.
Using \eqref{q:Syr} and \eqref{q:Syv} we can then bound
\begin{equation}\label{q:bdSy}
   |\Sy^*(t)|_0^{} \le c\,|\Rr^*|_s^{}
	\le c\,(|\qa(0)|^2 + \|f\|^2)\,\exp\bigl(-\sigma/\eps^{1/4}\bigr),
\end{equation}
also valid for $t\in[0,T_0]$.

The last term in \eqref{q:mudvhdt} can be bounded as
\begin{equation}\label{q:bdFy1}\begin{aligned}
   |\Fy^>|_0^{} &\le c\,|\gb_3W^<|\, |\gb_3W^>|
		\le c\,|\gb_3W|\, |\gb_3W^>|\\
		&\le c\,\exp(-\sigma/\eps^{1/4})\,
			|\gb_3W|\,\|\lapl3W\|_\sigma^{}\,.
\end{aligned}\end{equation}
We now take $C_1=2|\gb_3 W_0|$ and $C_\sigma=2\|\lapl3 W_0\|_\sigma^{}$
in Lemma~\ref{t:gev}.
Since $|\gb_3 W_0|\le c\,\|\lapl3 W_0\|_\sigma^{}$ by our hypothesis
(i.e.\ $\sigma>0$ and the rhs is finite), we have
\begin{equation}
   |\gb_3 W(t)| \le 2 |\gb_3 W_0|
   \qquad\textrm{and}\qquad
   \|\lapl3 W(t)\|_\sigma^{} \le 2 \|\lapl3 W_0\|_\sigma^{}
\end{equation}
for $t\in[0,T_\sigma']$ where
$T_\sigma'(\|\lapl3 W_0\|_\sigma^{},\|f\|;\sigma,\mu,\Dom)=\min(T_1,T_\sigma)$.
Thus \eqref{q:bdFy1} becomes
\begin{equation}\label{q:bdFy}
  |\Fy^>(t)|_0^{} \le c\,\exp(-\sigma/\eps^{1/4})\,\|\lapl3 W_0\|_\sigma^2
\end{equation}
for $t\in[0,T_\sigma']$.

Finally, writing $y=|\hat\vb|^2+|\hat\rho|^2$, we rewrite \eqref{q:mudvhdt} as
\begin{equation}\label{q:dydt}
   \ddt{y} \le M y + N
\end{equation}
where $M$ is bounded independently of $\eps$,
\begin{equation}\label{q:bdM}
   M(t) = \frac{c}\mu\, (|\qa(t)|_s^2 + \|f\|^2)
	\le \frac{2c}\mu\,(|\qa(0)|_s^2 + \|f\|^2)
	\le \frac{c'}\mu (\|\lapl3 q_0\|_\sigma^2 + \|f\|^2)
\end{equation}
for $t\in[0,T_0]$.
As for $N$, \eqref{q:bdSy} and \eqref{q:bdFy} give us
\begin{equation}\begin{aligned}
   N(t) &= \frac{c}\mu \bigl(|\Sy^*(t)|^2 + \Fy^>(t)|^2\bigr)\\
	&\le \cnst9(\mu)\, \bigl[ (|\gb_3^3 W_0|+\|f\|)^2
		+ \|\lapl3 W_0\|_\sigma^2 \bigr] \exp(-2\sigma/\eps^{1/4}),
\end{aligned}\end{equation}
valid for $t\in[0,T_\sigma'']$ where
$T_\sigma''(\|\gb_3^3 W_0\|_\sigma^{},\|f\|;\sigma,\mu,\Dom)=\min(T_0,T_\sigma')$.

Since $W^>$ and $\hat W$ are $L^2(\Dom)$-orthogonal, the hypothesis
\eqref{q:idhypo} implies that
\begin{equation}\label{q:Worth}
   |\hat W_0|^2 + |W^>_0|^2 \le C_{\rm id}\exp(-2\sigma/\eps^{1/4}).
\end{equation}
Lemma~\ref{t:gev} and \eqref{q:whisig} then imply that there is a
$T_3(\|\gb_3^3 W_0\|_\sigma^{},\|f\|;\sigma,\mu,\Dom)$ such that for
$t\in[0,T_3]$,
\begin{equation}
   |W^>(t)|^2 \le 2 C_{\rm id}\exp(-2\sigma/\eps^{1/4}).
\end{equation}
Integrating \eqref{q:dydt} and taking into account
\eqref{q:bdM}--\eqref{q:Worth} give us
\begin{equation}
   |\hat W(t)|^2 \le 2 C_{\rm id}\exp(-2\sigma/\eps^{1/4})
\end{equation}
for $t\in[0,T_4]$ with
$T_4=T_4(\|\gb_3^3 W_0\|_\sigma^{},\|f\|;\sigma,\mu,\Dom)$.
It thus follows that
\begin{equation}
   |W^>(t)|^2 + |\hat W(t)|^2 \le 4 C_{\rm id}\exp(-2\sigma/\eps^{1/4}),
\end{equation}
for $t\in[0,\min(T_3,T_4)]$, which is precisely \eqref{q:mainbd}.


\nocite{petcu:3dpe}
\nocite{allen:93}
\nocite{lions-temam-wang:92b}
\nocite{kobelkov:06}
\nocite{ju:3dpe}


\begin{thebibliography}{10}

\bibitem{allen:93}
{\sc J.~S. Allen}, {\em Iterated geostrophic intermediate models}, J. Phys.
  Ocean., 23 (1993), pp.~2447--2461.

\bibitem{bartello:95}
{\sc P.~Bartello}, {\em Geostrophic adjustment and inverse cascades in rotating
  stratified turbulence}, J. Atmos. Sci., 52 (1995), pp.~4410--4428.

\bibitem{bourgeois-beale:94}
{\sc A.~J. Bourgeois and J.~T. Beale}, {\em Validity of the quasigeostrophic
  model for large-scale flow in the atmosphere and ocean}, SIAM J. Math. Anal.,
  25 (1994), pp.~1023--1068.

\bibitem{cao-titi:u-3dpe}
{\sc C.~Cao and E.~S. Titi}, {\em Global well-posedness of the
  three-dimensional viscous primitive equations of large scale ocean and
  atmosphere dynamics}.
\newblock arXiv:math.AP/0503028v1, 2005.

\bibitem{cotter:th}
{\sc C.~J. Cotter}, {\em Model reduction for shallow water dynamics: balance,
  adiabatic invariance and subgrid modelling}, PhD thesis, Imperial College,
  London, 2004.

\bibitem{ju:3dpe}
{\sc N.~Ju}, {\em The global attractor for the solutions to the {3D} viscous
  {P}rimitive {E}quations}, Discrete Cont. Dyn. Systems, Ser.~A,  (2006), to
  appear.

\bibitem{kobelkov:06}
{\sc G.~M. Kobelkov}, {\em Existence of a solution ``in whole'' for the
  large-scale ocean dynamics equations}, C. R. Acad. Sc. Paris, S{\'e}r.~A,
  (2006), to appear.

\bibitem{lions-temam-wang:92b}
{\sc J.-L. Lions, R.~Temam, and S.~Wang}, {\em On the equations of the
  large-scale ocean}, Nonlinearity, 5 (1992), pp.~1007--1053.

\bibitem{lochak-meunier:mpavg}
{\sc P.~Lochak and C.~Meunier}, {\em Multiphase averaging for classical
  systems}, Springer-Verlag, 1988.

\bibitem{lorenz:80}
{\sc E.~N. Lorenz}, {\em Attractor sets and quasi-geostrophic equilibrium}, J.
  Atmos. Sci., 37 (1980), pp.~1685--1699.

\bibitem{matthies:01}
{\sc K.~Matthies}, {\em Time-averaging under fast periodic forcing of parabolic
  partial differential equations: exponential estimates}, J. Diff. Eq., 174
  (2001), pp.~133--180.

\bibitem{petcu:3dpe}
{\sc M.~Petcu}, {\em On the three dimensional primitive equations}, Adv. Diff.
  Eq.,  (2006), to appear.

\bibitem{petcu-temam-dw:rgpe}
{\sc M.~Petcu, R.~Temam, and D.~Wirosoetisno}, {\em Renormalization group
  method applied to the primitive equation}, J. Diff. Eq., 208 (2004),
  pp.~215--257.

\bibitem{petcu-dw:gev3}
{\sc M.~Petcu and D.~Wirosoetisno}, {\em Sobolev and {G}evrey regularity
  results for the primitive equations in 3 space dimensions}, Applic. Anal., 84
  (2005), pp.~769--788.

\bibitem{temam:iddsmp}
{\sc R.~Temam}, {\em Infinite-dimensional dynamical systems in mechanics and
  physics, 2nd ed.}, Springer-Verlag, 1997.

\bibitem{temam-dw:lbal}
{\sc R.~Temam and D.~Wirosoetisno}, {\em Long-time balance for the primitive
  equations of the ocean}.
\newblock Article in preparation, 2006.

\bibitem{jv-yavneh:04}
{\sc J.~Vanneste and I.~Yavneh}, {\em Exponentially small inertia--gravity
  waves and the breakdown of quasi-geostrophic balance}, J. Atmos. Sci., 61
  (2004), pp.~211--223.

\bibitem{wbsv:95}
{\sc T.~Warn, O.~Bokhove, T.~G. Shepherd, and G.~K. Vallis}, {\em Rossby number
  expansions, slaving principles, and balance dynamics}, Quart. J. Roy. Met.
  Soc., 121 (1995), pp.~723--739.

\bibitem{dw:abal}
{\sc D.~Wirosoetisno}, {\em Exponentially accurate balanced dynamics}, Adv.
  Diff. Eq., 9 (2004), pp.~177--196.

\bibitem{dw-tgs-temam:wwrg}
{\sc D.~Wirosoetisno, T.~G. Shepherd, and R.~M. Temam}, {\em Free gravity waves
  and balance dynamics}, J. Atmos. Sci., 59 (2002), pp.~3382--3398.

\end{thebibliography}

\end{document}